\documentclass[11pt,letterpaper,reqno]{amsart}

\usepackage[T1]{fontenc}
\usepackage[utf8]{inputenc} 
\usepackage{lmodern,microtype}
\usepackage[margin=1in]{geometry}

\usepackage{amsmath,amssymb,amsthm,mathtools}
\usepackage{esint,stmaryrd}
\usepackage{graphicx,overpic}
\usepackage{epsfig} 
\usepackage{caption}
\usepackage{xcolor,enumitem}
\usepackage[normalem]{ulem}
\usepackage[colorlinks=true,allcolors=blue!45!black]{hyperref}

\newcommand{\loc}{\mathrm{loc}} 
\numberwithin{equation}{section}
\numberwithin{figure}{section}

\theoremstyle{plain}
\newtheorem{theorem}{Theorem}[section]
\newtheorem{lemma}[theorem]{Lemma}

\newtheorem{corollary}[theorem]{Corollary}
\newtheorem*{theorem*}{Theorem}
\newtheorem*{corollary*}{Corollary}

\theoremstyle{definition}

\newtheorem{remark}[theorem]{Remark}

\newtheorem*{notation*}{Notation}

\allowdisplaybreaks
\setlist[enumerate]{itemsep=3pt,topsep=5pt}
\newcommand{\R}{\mathbb{R}}

\newcommand{\one}{\mathbf{1}} 

\newcommand{\Om}{\Omega}

\DeclareMathOperator{\diver}{div}

\newcommand{\Lip}{\operatorname{Lip}}

\newcommand{\restr}[2]{%
  \left.\kern-\nulldelimiterspace #1\right|_{#2}%
}

\title{On the regularity of stationary solutions to free boundary problems}

\author{Daniel Restrepo}
\address{127 Vincent Hall 206 Church St. SE Minneapolis, MN 55455, United States of America}
\email{drestrep@umn.edu}

\begin{document}
\maketitle

\begin{abstract}
In this note, we improve the regularity of H\"older continuous stationary solutions with respect to inner variations which are harmonic outside their zero set. We obtain a $\sqrt{\log}$-Lipschitz estimate for a class of solutions allowing a generalized phase variable and including the two-phase problem. In the case of the one-phase problem, we can improve to Lipschitz regularity. In particular, we show that the initial Lipschitz assumption in the class of variational solutions considered by Kriventsov and Weiss can be replaced by any positive H\"older exponent.
\end{abstract}

\section{Introduction}
Given an open set $\Omega\subset\R^n$, we consider stationary points of the energy
\begin{equation}\label{eq energy}
 J(u)=\int_\Omega (|\nabla u|^2+\chi)\,dx,
\end{equation}
with respect to inner variations (or domain variations) which are also harmonic in the complement of their zero set. Here $\chi=\chi(u)$ is bounded and measurable. Functions satisfying these conditions arise as critical points of elliptic free boundary problems, such as the one-phase problem \cite{AC}, where $\chi=\one_{\{u>0\}}$, or the two-phase problem, where $\chi=\lambda_1\one_{\{u>0\}}+\lambda_2\one_{\{u<0\}}$ with $\lambda_i>0$ for $i=1,2$ \cite{ACF}. Related stationarity conditions also arise naturally in optimal partition and segregation problems \cite{terraciniverziniconti,CaffarelliLin} and in semilinear problems with a topological constraint \cite{NRS}.

In the one-phase problem, considering critical points with respect to domain variations instead of outer variations has the advantage of admitting solutions which are neither minimizers nor critical points with respect to outer variations, such as $v(x)=|x_1|$. More importantly, the compactness theory developed by Kriventsov and Weiss \cite{KW} provides stability under limits, in contrast with the notion of weak solution introduced in \cite{AC}.

Regularity for weak notions of solution to the one-phase problem is subtle. Lewis and Vogel \cite{LV91,LV} constructed examples of distributional solutions which are H\"older continuous but not Lipschitz: their positivity sets have locally finite perimeter, and the solutions satisfy the distributional formulation of the Bernoulli boundary condition underlying the notion of weak solution introduced by Alt and Caffarelli \cite{AC,DPR}. In \cite{AC}, Alt and Caffarelli develop a regularity theory for weak solutions by imposing regularity and nondegeneracy assumptions from the beginning. The examples of Lewis and Vogel show that such hypotheses are necessary since non-Lipschitz distributional solutions exist. Similarly, in their notion of variational solution, Kriventsov and Weiss assume Lipschitz regularity a priori and conjecture that stationarity with respect to inner variations is not enough to guarantee Lipschitz continuity. Our main result shows that every H\"older continuous function satisfying the remaining conditions in their definition of variational solution is necessarily Lipschitz.

\begin{theorem}\label{thm main}
Let $\Om\subset\R^n$ be open, $n\geq2$, and $\alpha\in(0,1)$. Suppose that $u:\Om\to[0,\infty)$ and $\chi:\Om\to[0,1]$ satisfy
\begin{enumerate}[label=(\roman*)]
\item $u\in H^1_\loc( \Om)\cap  C^\alpha_\loc(\Om)$ and $u$ is harmonic in $\{u>0\}$;
\item $\chi$ is measurable and $\one_{\{u>0\}}\leq\chi$ almost everywhere;
\item for every $T\in\Lip_c(\Om;\R^n)$,
\begin{equation}\label{eq stationarity}
 \int_\Om\left[(|\nabla u|^2+\chi)\diver T-2\nabla u\cdot DT\nabla u\right]dx=0.
\end{equation}
\end{enumerate}
Then $u$ is locally Lipschitz. More precisely, if $B_{2R}(x_0)\subset \subset\Om$, then
\begin{equation}\label{eq mainbound}
 \|\nabla u\|_{L^\infty(B_R(x_0))}
 \leq C\left(1+\frac{\|u\|_{L^\infty(B_{2R}(x_0))}}R\right).
\end{equation}
\end{theorem}

Some remarks are in order. First, the initial H\"older regularity plays a qualitative role, providing the compactness needed in the proof. Neither $\alpha$ nor the initial H\"older seminorm appears in \eqref{eq mainbound}. Second, the theorem also allows $\chi$ to be a generalized phase variable: it need not be constant on the level sets of $u$, and, apart from the condition $\one_{\{u>0\}}\leq\chi$, the two are decoupled a priori; see \cite{KW} for further discussion. In our theorem, $\chi$ may take any value in $[0,1]$, rather than only the values $0$ and $1$ allowed in the definition of Kriventsov and Weiss. This larger class is natural when phases are only known to converge weakly-*, since such limits retain boundedness but need not retain a binary range. Third, the improvement to $\sqrt{\log}$-Lipschitz regularity holds for sign-changing $u$ and arbitrary $\chi\in L^\infty_\loc(\Om)$; see Theorem~\ref{thm anyholder}. The final step to Lipschitz regularity uses both nonnegativity and the comparison between $\chi$ and the positivity phase.

Optimal Lipschitz regularity for the two-phase problem seems more delicate, and a suitable notion of variational solution with a relaxed phase deserves more careful analysis.

Theorem~\ref{thm main} shows directly that the initial Lipschitz assumption for variational solutions in \cite[Section 3.2]{KW} can be lowered to $\alpha$-H\"older regularity for any $\alpha\in(0,1)$. Moreover, \eqref{eq mainbound} provides a Lipschitz bound in terms of an $L^\infty$ bound. In recent work, De Philippis--Hirsch--Nahon \cite[Corollary 7]{DHN} prove local Lipschitz regularity in dimension two for $H^1$ functions satisfying \eqref{eq stationarity} with $\chi=\one_{\{u\neq0\}}$, without assuming harmonicity in the nonzero set; see Section~\ref{sec harmonic}.

Finally, Theorem~\ref{thm main} shows that the Lewis--Vogel examples cannot satisfy the inner variation condition \eqref{eq stationarity} with $\chi=\one_{\{u>0\}}$; otherwise, the would be Lipschitz.

\subsection{On the sharpness of the assumptions}\label{sec harmonic}
We discuss first the assumption $\Delta u=0$ in $\{u>0\}$. As proved in \cite[Proposition 3.2]{KW}, the stationarity identity \eqref{eq stationarity} implies harmonicity in the positivity set if $u\in C^2_\loc(\{u>0\})$. Indeed, for $T\in C_c^\infty(\{u>0\};\R^n)$,
\begin{align}
 \diver\bigl(|\nabla u|^2T-2(T\cdot\nabla u)\nabla u\bigr)
 &=|\nabla u|^2\diver T-2\nabla u\cdot DT\nabla u \notag\\
 &\quad-2(T\cdot\nabla u)\Delta u.
 \label{eq dividentity}
\end{align}
Since $\chi=1$ on the support of $T$, integration and \eqref{eq stationarity} give
$\int(T\cdot\nabla u)\Delta u=0$. Thus $\Delta u=0$ wherever
$\nabla u\neq0$. By continuity of $\Delta u$, this also holds on the
closure of that set; on the interior of $\{\nabla u=0\}$, the function
is locally constant.

One can weaken $C^2_\loc(\{u>0\})$ to $C^1_\loc(\{u>0\})$ when $\Delta u$ is a Radon measure in the positivity set. To see this, write $\mu=\Delta u$ and mollify. Then $\nabla u_\varepsilon\to\nabla u$ uniformly on compact subsets of
$\{u>0\}$, while $\Delta u_\varepsilon\rightharpoonup^*\mu$ as Radon
measures. Integrating \eqref{eq dividentity}, passing to the limit, and using \eqref{eq stationarity}, we obtain
\begin{equation*}
 \int_{\{u>0\}}(T\cdot\nabla u)\,d\mu=0
\end{equation*}
for every test $T$. Therefore, $u$ is harmonic on
$\{u>0\}\cap\{|\nabla u|>0\}$. Moreover, one can verify in the viscosity sense that $u$ solves $|\nabla u|^2\Delta u=0$ in $\{u>0\}$. By a result of Imbert and Silvestre \cite[Lemma 6]{IS}, it follows that $\Delta u=0$ in $\{u>0\}$.

Lipschitz regularity alone is insufficient to show harmonicity in the positivity set. Since \eqref{eq stationarity} is quadratic in the gradient, one can produce solutions of \eqref{eq stationarity} with $\chi=1$ by taking $u=|\phi|+1$ with $\phi$ harmonic. In particular, $u(x)=1+|x_1|$ is a Lipschitz solution of \eqref{eq stationarity}.
Without extra regularity in the positivity set, stationarity does not distinguish the zero level of $u$. In dimension two, if $u$ satisfies \eqref{eq stationarity} with $\chi$ constant, then $u$ is automatically locally Lipschitz. If, furthermore, $\Delta u$ is a signed Radon measure, its support is locally constrained to level sets of $u$; see \cite{Rodiac}.

When $n\geq3$, it is not known whether regularity of $u$ can be obtained directly from \eqref{eq stationarity}, even when $\chi$ is constant. If one assumes in addition that $u\Delta u=0$ locally, in the distributional form \eqref{eq outer}, one recovers monotonicity of the Almgren frequency $N(z,r)$ when $\chi=0$; see \eqref{eq frequency0}. In the presence of the phase $\chi$, \eqref{eq outer} and \eqref{eq stationarity} together give the monotonicity of the adjusted frequency \eqref{eq newmon} while its value is below one. This monotonicity formulae makes equivalent initial H\"older continuity and (local) uniform lower frequency bounds. A lower frequency bound can also be obtained from geometric information: for instance, when $\{u>0\}\cap B_r(z)$ has at least two connected components of comparable size for all sufficiently small $r$. This observation is exploited in \cite{NRS}, where $\chi=f(u)$ with $f$ Lipschitz and $f(0)=0$, to derive lower bounds for $N(z,r)$ from a notion of minimality.

For the Dirichlet case ($\chi=0$), one can also use monotonicity and a covering argument to bound the Hausdorff dimension of points with small frequency. More precisely, for $\beta>0$ one has
\begin{equation*}
 \dim_{\mathcal H}
 \bigl\{z\in\Om:u(z)=0,\ \lim_{r\downarrow0}N(z,r)<\beta\bigr\}
 \leq n-2+2\beta,
\end{equation*}
which gives
\begin{equation*}
 \dim_{\mathcal H}
 \bigl\{z\in\Om:u(z)=0,\ \lim_{r\downarrow0}N(z,r)=0\bigr\}
 \leq n-2.
\end{equation*}
Sets of finite $\mathcal H^{n-2}$ measure have zero Dirichlet capacity \cite[Section 4.7.2]{EvansGariepy}, but a dimension bound at $n-2$ alone does not guarantee zero capacity. A refinement showing that the small-frequency set has zero capacity would make it invisible to the Dirichlet energy and could allow one to remove it and deduce H\"older continuity. In this regard, an example of a variational solution harmonic outside its zero set but lacking H\"older continuity would lie at the borderline between sets of zero capacity and sets of dimension $n-2$ with positive capacity.

\subsection{Proof strategy}\label{sec proofstrategy}

The first part follows the compactness and frequency arguments of \cite[Section 2]{NRS}. For the Dirichlet problem, $\chi=0$, Almgren monotonicity and the initial H\"older regularity give compactness of height rescalings. The limits are homogeneous, and dimension reduction gives a lower frequency bound of one. We include the argument in Section~\ref{subsec frequency}.

For a nonzero phase, we adapt this strategy using an adjusted frequency;  see Lemma \ref{thm monF}. This allows us to obtain a $\sqrt{\log}$-Lipschitz estimate for a larger class of solutions, including sign-changing $u$ and locally bounded $\chi$. This is the content of Theorem \ref{thm anyholder}. The formula is related to the frequency calculations in \cite{KW}, but is used here before Lipschitz regularity is known.

Passing from $\sqrt{\log}$-Lipschitz regularity to Lipschitz regularity is subtler because the adjusted frequency need not be monotone when its value reaches one. In the one-phase case, we first prove that the linear rescalings have bounded height (normalized $L^2$ norm) at each fixed free boundary point. This gives the existence of the generalized phase density, namely the limit of the averages of $\chi$. At points of maximal generalized density, the height of Lipscthiz rescalings is monotone. At the remaining points, we show that $u$ satisfies a viscosity supersolution condition for the one-phase problem. These two observations give the uniform Lipschitz estimate.

\subsection{Organization of the paper}

Section~\ref{sec holder} is devoted to the proof of Theorem~\ref{thm anyholder}, which provides a $\sqrt{\log}$-Lipschitz improvement for solutions with a general locally bounded phase $\chi$. In Section~\ref{sec lipschitz}, we prove Theorem~\ref{thm main}.

\subsection{Notation and conventions}

Throughout the paper, $C$ denotes a generic positive constant depending only on $n$, whose value may increase from line to line; and whose dependencies will be clearly established in the statements of
the corresponding propositions. We use $C^\alpha$ for H\"older spaces with $0<\alpha<1$ and $\Lip$ for Lipschitz spaces; $C^1$ denotes continuously differentiable functions. For an open set $U$, our H\"older seminorm is
\begin{equation*}
 [u]_{C^\alpha(U)}
 =\sup_{\substack{x,y\in U\\x\neq y}}
 \frac{|u(x)-u(y)|}{|x-y|^\alpha}.
\end{equation*}

\medskip

\noindent
\textbf{Acknowledgments.} The author would like to thank to Max Engelstein and Dennis Kriventsov for some insightful comments on a preliminary version of this work. The author has been supported by the AMS-Simons travel grant.

\medskip

\noindent
\textbf{Statement of AI usage.} No LLM's were used in the generation of the mathematical content of the paper. AI tools were used for improvements of phrase and construction.

\section{Improvement of H\"older continuity}\label{sec holder}

As announced in the introduction, the results of this section hold for a class slightly larger class of phases and for sign-changing functions.  
\begin{theorem}\label{thm anyholder}
Let $\Om\subset\R^n$ be open, $n\geq2$, and let $u\in H^1_\loc(\Om)\cap C^\alpha_\loc(\Om)$ for some $\alpha\in(0,1)$. Assume that $u$ is harmonic in $\{u\neq0\}$, $\chi\in L^\infty_\loc(\Om)$, and
\begin{equation}\label{eq lambda}
 \int_\Om\left[(|\nabla u|^2+\chi)\diver T
 -2\nabla u\cdot DT\nabla u\right]dx=0
 \quad\mbox{for all }T\in\Lip_c(\Om;\R^n).
\end{equation}
Then $u\in C^\gamma_\loc(\Om)$ for every $\gamma\in(0,1)$. More precisely, if $B_{2R}(x_0)\subset \subset\Om$, $x,y\in B_R(x_0)$, and $0<|x-y|\leq R/2$, then
\begin{equation}\label{eq anyHolderbound}
\begin{split}
 |u(x)-u(y)|\leq C|x-y|\bigg(&\frac{\|u\|_{L^\infty(B_{2R}(x_0))}}R\\
 &+\|\chi\|_{L^\infty(B_{2R}(x_0))}^{1/2}
 \sqrt{1+\log\frac{R}{|x-y|}}\bigg).
\end{split}
\end{equation}
In particular, $u$ is locally $\sqrt{\log}$-Lipschitz.
\end{theorem}

The assumption $\chi\geq 0$ is not restrictive. Indeed, adding a constant to $\chi$ does not change \eqref{eq lambda}, since $\int\diver T=0$. On any fixed relatively compact ball, we may therefore replace $\chi$ by $\chi+\|\chi\|_{L^\infty}$.

\subsection{Basic identities and monotonicity formulae}\label{subsec monotonicity}

Throughout the rest of this section, we assume that $u$ and $\chi$ satisfy the hypotheses of Theorem~\ref{thm anyholder}. We write $u=u^++u^-$, where $u^+=\max\{u,0\}$ and $u^-=\min\{u,0\}$. We observe that both $u^+$ and $-u^-$ are subharmonic, since they are nonnegative and harmonic in their positivity sets. Indeed, their harmonic extensions on a ball are nonnegative, and the maximum principle on each positivity component gives the comparison. Consequently, $\Delta u^+=\mu^+\geq0$ and $\Delta(-u^-)=\mu^-\geq0$ are Radon measures supported on $\partial\{u>0\}$ and $\partial\{u<0\}$, respectively. In particular, $\Delta u=\mu^+-\mu^-$ is a signed Radon measure supported on $\partial\{u\neq0\}$.

We also have $u\Delta u=0$ in the sense of distributions in $\Om$, that is,
\begin{equation}\label{eq outer}
 \int_\Om\eta|\nabla u|^2
 =-\int_\Om u\nabla u\cdot\nabla\eta
 \quad\mbox{for all }\eta\in C_c^\infty(\Om).
\end{equation}
To see this, test harmonicity in $\{u>0\}$ with $\eta(u-\varepsilon)_+$ and in $\{u<0\}$ with the analogous truncation of $-u$, then let $\varepsilon\downarrow0$. Here we use that $\nabla u=0$ almost everywhere on $\{u=0\}$. It follows as well that $\Delta(u^2)=2|\nabla u|^2$ in the sense of distributions. Taking cutoffs converging to $\one_{B_r(z)}$ in \eqref{eq outer}, we obtain, for almost every admissible radius,
\begin{equation}\label{eq testdir}
 \int_{B_r(z)}|\nabla u|^2
 =\int_{\partial B_r(z)}u\,\partial_\nu u,
\end{equation}
where $\nu$ is the outer unit normal to the sphere.

Similarly, testing \eqref{eq lambda} with radial vector fields of the form $T(x)=(x-z)\eta_\varepsilon(x)$, where $\eta_\varepsilon\to\one_{B_r(z)}$, yields the Pohozaev identity
\begin{equation}\label{eq pohozaev}
 (n-2)\int_{B_r(z)}|\nabla u|^2+n\int_{B_r(z)}\chi
 =r\int_{\partial B_r(z)}
 \bigl(|\nabla u|^2+\chi-2(\partial_\nu u)^2\bigr)
\end{equation}
for almost every $r>0$ such that $B_r(z)\subset \subset\Om$; compare \cite[Proposition 3.4]{KW}.

Combining these identities we show that the Weiss-type energies $W_\beta$, addapted for $\beta$-H\"older rescalings, are monotone. We reserve the notation $W(z,r)=W_1(z,r)$ for the standard Lipschitz Weiss energy; see also \cite{Weiss99,KW}.

\begin{lemma}[Weiss monotonicity]\label{lemma weissbeta}
For $0<\beta\leq1$, set
\begin{equation}\label{eq weissbeta}
 W_\beta(z,r)
 =r^{2-n-2\beta}\int_{B_r(z)}(|\nabla u|^2+\chi)
 -\beta r^{1-n-2\beta}\int_{\partial B_r(z)}u^2.
\end{equation}
Then, for almost every admissible radius,
\begin{equation}\label{eq weissbeta derivative}
\begin{split}
 \frac{d}{dr}W_\beta(z,r)
 ={}&2r^{2-n-2\beta}
 \int_{\partial B_r(z)}
 \left(\partial_\nu u-\frac{\beta u}{r}\right)^2\\
 &+2(1-\beta)r^{1-n-2\beta}\int_{B_r(z)}\chi
 \geq0.
\end{split}
\end{equation}
If $u(z)=0$ and $\beta<\alpha$, then $W_\beta(z,0+)=0$, and hence
\begin{equation}\label{eq Fbeta}
 \int_{B_r(z)}(|\nabla u|^2+\chi)
 \geq\frac{\beta}{r}\int_{\partial B_r(z)}u^2.
\end{equation}
\end{lemma}

\begin{proof}
Write $E(r)=\int_{B_r(z)}(|\nabla u|^2+\chi)$ and
$S(r)=\int_{\partial B_r(z)}u^2$. Differentiating the first term in
$W_\beta$ and using \eqref{eq pohozaev}, we find
\begin{equation}\label{eq first term}
\begin{split}
 \frac{d}{dr}\bigl(r^{2-n-2\beta}E(r)\bigr)
 =2r^{1-n-2\beta}\bigg[
 &(1-\beta)\int_{B_r(z)}\chi
 -\beta\int_{B_r(z)}|\nabla u|^2\\
 &+r\int_{\partial B_r(z)}(\partial_\nu u)^2
 \bigg].
\end{split}
\end{equation}
The other term can be written in terms of the
$\beta$-H\"older rescalings $u_r(x)=u(z+rx)/r^\beta$ as
\begin{equation*}
 r^{1-n-2\beta}S(r)=\int_{\partial B_1}u_r^2.
\end{equation*}
Differentiating gives
\begin{equation}\label{eq second term}
 \frac{d}{dr}\bigl(r^{1-n-2\beta}S(r)\bigr)
 =2r^{1-n-2\beta}
 \int_{\partial B_r(z)}
 u\left(\partial_\nu u-\frac{\beta u}{r}\right).
\end{equation}
We now subtract $\beta$ times \eqref{eq second term} from
\eqref{eq first term}. By \eqref{eq testdir}, the integral of
$|\nabla u|^2$ in \eqref{eq first term} becomes an integral of
$u\partial_\nu u$ over the sphere. The resulting boundary terms
complete the square in \eqref{eq weissbeta derivative}.

If $u(z)=0$, H\"older continuity gives
$\int_{\partial B_r(z)}u^2=O(r^{n-1+2\alpha})$. Testing
\eqref{eq outer} with $\eta^2$, where $\eta$ is a suitable cutoff,
yields the Caccioppoli inequality
\begin{equation*}
 \int_{B_r(z)}|\nabla u|^2
 \leq\frac{C}{r^2}\int_{B_{2r}(z)}u^2
 =O(r^{n-2+2\alpha}).
\end{equation*}
Since $\chi$ is locally bounded,
$\int_{B_r(z)}\chi=O(r^n)$. Thus both terms in
$W_\beta(z,r)$ tend to zero when $\beta<\alpha$, and
\eqref{eq Fbeta} follows from monotonicity.
\end{proof}

For completeness, we next show that Almgren's frequency is monotone when $\chi=0$. The corresponding calculation in \cite[Lemma 2.9]{NRS} is somewhat more involved because of the potential term. 

Given $z\in\Om$, an admissible radius $r$, and $v\in H^1_\loc(\Om)$, define the scaled Dirichlet energy
\begin{equation}\label{eq Dir}
 D_v(z,r)=r^{2-n}\int_{B_r(z)}|\nabla v|^2,
\end{equation}
the height
\begin{equation}\label{eq Height}
 H_v(z,r)=r^{1-n}\int_{\partial B_r(z)}v^2,
\end{equation}
and, whenever $H_v(z,r)>0$, the frequency
\begin{equation}\label{eq frequency0}
 N_v(z,r)=\frac{D_v(z,r)}{H_v(z,r)}
 =\frac{r\int_{B_r(z)}|\nabla v|^2}
        {\int_{\partial B_r(z)}v^2}.
\end{equation}
When the function is fixed, we omit its subscript and write
$D(z,r)$, $H(z,r)$, and $N(z,r)$.

\begin{remark}\label{remark frequencylog}
For a function $v\in H^1_\loc(\Om)$ satisfying $v\Delta v=0$ in the distributional sense, differentiation of the height and the analogue of \eqref{eq testdir} give
\begin{equation}\label{eq derivativeheight}
 \frac{d}{dr}H_v(z,r)=\frac{2}{r}D_v(z,r).
\end{equation}
Consequently, wherever the height is positive,
\begin{equation}\label{eq logheightder}
 N_v(z,r)=\frac{r}{2}\frac{d}{dr}\log H_v(z,r).
\end{equation}
\end{remark}

\begin{lemma}[Almgren monotonicity]\label{lemma almgren}
If $\chi=0$, then $N(z,r)$ is nondecreasing wherever $H(z,r)>0$. More precisely,
\begin{equation}\label{eq freqder}
 \frac{d}{dr}N(z,r)
 =\frac{2r^{2-n}}{H(z,r)}
 \int_{\partial B_r(z)}
 \left(\partial_\nu u-\frac{N(z,r)}{r}u\right)^2
 \geq0.
\end{equation}
If $N(z,r)$ is constant on $(r_1,r_2)$, then $u$ is homogeneous about $z$ on $B_{r_2}(z)\setminus\overline{B_{r_1}(z)}$, of degree equal to that constant.
\end{lemma}

\begin{proof}
Since the frequency is invariant under rescaling, it suffices, after translation and scaling, to compute the derivative at $z=0$ and $r=1$. Differentiating $N=D/H$ gives
\begin{equation}\label{eq derfrequ}
 N'(0,1)
 =\frac{D'(0,1)H(0,1)-H'(0,1)D(0,1)}{H(0,1)^2}.
\end{equation}
Equation \eqref{eq first term} with $\chi=0$ and $\beta=0$ gives
\begin{equation}\label{eq derdir}
 D'(0,1)=2\int_{\partial B_1}(\partial_\nu u)^2,
\end{equation}
while \eqref{eq derivativeheight} and \eqref{eq testdir} give
\begin{equation}\label{eq derheight}
 H'(0,1)=2D(0,1)
 =2\int_{\partial B_1}u\,\partial_\nu u.
\end{equation}
Thus the numerator in \eqref{eq derfrequ} is twice the
Cauchy--Schwarz remainder
\begin{equation*}
 \left(\int_{\partial B_1}u^2\right)
 \left(\int_{\partial B_1}(\partial_\nu u)^2\right)
 -\left(\int_{\partial B_1}u\,\partial_\nu u\right)^2.
\end{equation*}
Writing $N(0,1)=D(0,1)/H(0,1)$, we can use \eqref{eq testdir} to write $D(0,1)=\int_{\partial B_1}u\,\partial_\nu u$ and express the remainder as
\begin{align}\notag
 N'(0,1)
 &=\frac{2}{H(0,1)}
   \left(\int_{\partial B_1}(\partial_\nu u)^2
         -\frac{D(0,1)^2}{H(0,1)}\right)\\ \label{eq freqcompu}
 &=\frac{2}{H(0,1)}
   \int_{\partial B_1}
   \bigl(\partial_\nu u-N(0,1)u\bigr)^2.
\end{align}
Rescaling gives \eqref{eq freqder}. Finally, rigidity follows from \eqref{eq freqder}.
\end{proof}

We finish this section with an adjusted frequency formula related to the one in \cite[Lemma 4.1]{KW}. In this case, however, monotonicity holds when the frequency is less than one, which is complementary to the approach in \cite{KW}, where the solutions considered have frequency at least one on the free boundary by the Lipschitz regularity assumption.

\begin{theorem}[Adjusted frequency]\label{thm monF}
Let $z\in\partial\{u\neq 0\}$. The adjusted frequency
\begin{equation}\label{eq newmon}
 F(z,r)
 =\frac{r\int_{B_r(z)}(|\nabla u|^2+\chi)}
        {\int_{\partial B_r(z)}u^2}
\end{equation}
satisfies
\begin{equation}\label{eq dermodfreq}
\begin{split}
 \frac{d}{dr}F(z,r)
 ={}&\frac{2r^{2-n}}{H(z,r)}
 \int_{\partial B_r(z)}
 \left(\partial_\nu u-\frac{N(z,r)}{r}u\right)^2\\
 &+\frac{2\int_{B_r(z)}\chi}
         {\int_{\partial B_r(z)}u^2}
 (1-N(z,r)).
\end{split}
\end{equation}
In particular, $r\mapsto\min\{F(z,r),1\}$ is nondecreasing.
\end{theorem}

\begin{proof}
With $E(r)$ as in the proof of Weiss monotonicity, we have
\begin{equation}\label{eq MF}
 F(z,r)=\frac{r^{2-n}E(r)}{H(z,r)}.
\end{equation}
The $\beta=0$ case of \eqref{eq first term} gives
\begin{equation}\label{eq firstder}
 \frac{d}{dr}\bigl(r^{2-n}E(r)\bigr)
 =2r^{2-n}\int_{\partial B_r(z)}(\partial_\nu u)^2
 +2r^{1-n}\int_{B_r(z)}\chi.
\end{equation}
Using \eqref{eq derivativeheight} to differentiate the quotient in
\eqref{eq MF}, we obtain
\begin{equation}\label{eq QR}
\begin{split}
 F'(z,r)
 ={}&\frac{2r^{2-n}}{H(z,r)}
       \int_{\partial B_r(z)}(\partial_\nu u)^2\\
 &+\frac{2\int_{B_r(z)}\chi}
         {\int_{\partial B_r(z)}u^2}
 -\frac{2}{r}F(z,r)N(z,r).
\end{split}
\end{equation}
Since
\begin{equation*}
 F(z,r)-N(z,r)
 =\frac{r\int_{B_r(z)}\chi}
        {\int_{\partial B_r(z)}u^2},
\end{equation*}
we have
\begin{equation}\label{eq inter}
 \frac{2\int_{B_r(z)}\chi}
      {\int_{\partial B_r(z)}u^2}
 -\frac{2}{r}F(z,r)N(z,r)
 =
 \frac{2\int_{B_r(z)}\chi}
      {\int_{\partial B_r(z)}u^2}
 (1-N(z,r))
 -\frac{2}{r}N(z,r)^2.
\end{equation}
The last term in \eqref{eq inter} combines with the first term on the right-hand side of \eqref{eq QR} to give the same square as in \eqref{eq freqcompu}. Indeed, \eqref{eq testdir} gives
\begin{equation*}
 \frac{2r^{2-n}}{H(z,r)}
 \int_{\partial B_r(z)}(\partial_\nu u)^2
 -\frac{2}{r}N(z,r)^2
 =
 \frac{2r^{2-n}}{H(z,r)}
 \int_{\partial B_r(z)}
 \left(\partial_\nu u-\frac{N(z,r)}{r}u\right)^2.
\end{equation*}
Together with \eqref{eq inter}, this proves \eqref{eq dermodfreq}.

Finally, $\chi\geq0$ implies $N(z,r)\leq F(z,r)$. Wherever
$F(z,r)<1$, we therefore have $N(z,r)<1$, so both terms in
\eqref{eq dermodfreq} are nonnegative. We conclude that
$r\mapsto\min\{F(z,r),1\}$ is nondecreasing.
\end{proof}

\subsection{Estimates and compactness}\label{subsec compactness}

In this subsection, we derive a regularity criterion from bounds on the height and then use it to prove strong compactness. The criterion will be used later both to improve the initial H\"older bound and, at the endpoint, to obtain Lipschitz regularity.

\begin{lemma}\label{lemma reg criterion}
Let $\gamma\in(0,1]$ and $K\geq0$. Assume that
\begin{equation}\label{eq heightbound}
 H(z,r)\leq Kr^{2\gamma}
\end{equation}
for every $z\in\partial\{u \neq 0\}$ and every $r>0$ such that $B_r(z)\subset\subset\Om$. Then $u\in C^\gamma_\loc(\Om)$ if $\gamma<1$, while $u\in\Lip_\loc(\Om)$ if $\gamma=1$. Moreover, whenever $B_{2R}(x_0)\subset\subset\Om$,
\begin{equation}\label{eq holderestimate}
 R^\gamma
 \sup_{\substack{x,y\in B_{R/2}(x_0)\\x\neq y}}
 \frac{|u(x)-u(y)|}{|x-y|^\gamma}
 \leq C\left(
 \sqrt K\,R^\gamma+
 R^{-n/2}\|u\|_{L^2(B_{2R}(x_0))}
 \right).
\end{equation}
\end{lemma}

\begin{proof}
Since $\Delta(u^2)=2|\nabla u|^2\geq0$, the function $u^2$ is subharmonic. Its mean-value inequality and the monotonicity of spherical averages give
\begin{equation}\label{eq heighttosup}
 \sup_{B_{r/2}(z)}|u|^2
 \leq C H(z,r)
 \leq CKr^{2\gamma}
\end{equation}
at every free boundary point $z$ for which $B_r(z)\subset\subset\Om$.

We first estimate $u$ and its gradient near the zero set. Let $x\in B_R(x_0)\cap\{u\neq0\}$, and set $d=\operatorname{dist}(x,\{u=0\})$. If $d<R/8$, a zero $z$ nearest to $x$ belongs to the boundary of the phase containing $x$. Applying \eqref{eq heighttosup} at $z$ with a radius comparable to $d$, and then applying the interior gradient estimate to the harmonic function $u$ in $B_d(x)$, we obtain
\begin{equation}\label{eq nearzero}
 |u(x)|\leq C\sqrt K\,d^\gamma,
 \qquad
 |\nabla u(x)|\leq C\sqrt K\,d^{\gamma-1}.
\end{equation}
On the other hand, where $d\geq R/8$, interior harmonic estimates and the mean-value inequality for $u^2$ give
\begin{equation}\label{eq awayzero}
 |\nabla u(x)|
 \leq CR^{-1-n/2}\|u\|_{L^2(B_{2R}(x_0))}.
\end{equation}

Now take $x,y\in B_{R/2}(x_0)$ and write $s=|x-y|$. If $s$ is comparable to $R$, the mean-value inequality for $u^2$ gives \eqref{eq holderestimate}. We may therefore assume that $s$ is small compared to $R$.

If either endpoint is within distance $2s$ of the zero set, both endpoints are within distance $3s$ of it. The first estimate in \eqref{eq nearzero}, applied in their respective phases, gives
\begin{equation*}
 |u(x)-u(y)|\leq C\sqrt K\,s^\gamma.
\end{equation*}
The same conclusion holds if an endpoint lies in the zero set. Otherwise, both endpoints are farther than $2s$ from the zero set, so the segment joining them stays at distance at least $s$ from it. We can integrate the gradient along this segment, using \eqref{eq nearzero} where its distance from the zero set is small and \eqref{eq awayzero} elsewhere. Since $\gamma\leq1$, this yields
\begin{equation*}
 |u(x)-u(y)|
 \leq C\sqrt K\,s^\gamma
 +CsR^{-1-n/2}\|u\|_{L^2(B_{2R}(x_0))}.
\end{equation*}
As $s\leq R$, this proves \eqref{eq holderestimate}. The constant is uniform for $\gamma\in(0,1]$, including the Lipschitz endpoint.
\end{proof}

As a corollary, we show that a lower bound for the adjusted frequency $F$ gives a H\"older estimate. The presence of the phase $\chi$ creates an error term that prevents us from getting Lipschitz regularity from this estimate.

\begin{corollary}\label{cor weisstoholder}
Let $\gamma\in(0,1)$. Assume that
\begin{equation}\label{eq remainderWeiss}
 \int_{B_r(z)}\bigl(|\nabla u|^2+\chi\bigr)
 \geq\frac{\gamma}{r}\int_{\partial B_r(z)}u^2
\end{equation}
for every $z\in\partial\{u \neq 0\}$ and every $r>0$ such that $B_r(z)\subset\subset\Om$. Then $u\in C^\gamma_\loc(\Om)$ and, whenever $B_{2R}(x_0)\subset\subset\Om$,
\begin{equation}\label{eq holderestimate2}
 R^\gamma[u]_{C^\gamma(B_{R/2}(x_0))}
 \leq C\left(
 R^{-n/2}\|u\|_{L^2(B_{2R}(x_0))}
 +\frac{R}{\sqrt{1-\gamma}}\,
 \|\chi\|_{L^\infty(B_{2R}(x_0))}^{1/2}
 \right).
\end{equation}
\end{corollary}

\begin{proof}
Write $M=\|\chi\|_{L^\infty(B_{2R}(x_0))}$. By \eqref{eq derivativeheight} and \eqref{eq remainderWeiss}, at every free boundary point $z$ and for admissible radii contained in $B_{2R}(x_0)$,
\begin{align*}
 \frac{d}{dr}\left(\frac{H(z,r)}{r^{2\gamma}}\right)
 &=\frac{2}{r^{2\gamma+1}}\bigl(D(z,r)-\gamma H(z,r)\bigr)\\
 &\geq-2r^{1-n-2\gamma}\int_{B_r(z)}\chi
 \geq-2|B_1|M r^{1-2\gamma}.
\end{align*}
Integrating from $s$ to $t$ gives
\begin{equation}\label{eq heightfromfrequency}
 \frac{H(z,s)}{s^{2\gamma}}
 \leq\frac{H(z,t)}{t^{2\gamma}}
 +\frac{|B_1|}{1-\gamma}M t^{2-2\gamma},
 \qquad 0<s<t.
\end{equation}

For the free boundary points needed to estimate $u$ in $B_{R/2}(x_0)$, take $t$ comparable to $R$. The mean-value inequality for $u^2$ gives $H(z,t)\leq CR^{-n}\|u\|_{L^2(B_{2R}(x_0))}^2$. This observation combined with \eqref{eq heightfromfrequency} provides the uniform height bound required to apply Lemma~\ref{lemma reg criterion}. Finally, this yields \eqref{eq holderestimate2}.
\end{proof}

\begin{lemma}[Strong compactness]\label{lemma compactness}
Let $(u_k,\chi_k)$ be a sequence of pairs satisfying the hypotheses of Theorem~\ref{thm anyholder} in a domain $D$, with a common exponent $\alpha\in(0,1)$. Assume that
\begin{equation*}
 \sup_k\|u_k\|_{H^1(D)}<\infty,
 \qquad
 \sup_k\|\chi_k\|_{L^\infty(D)}<\infty.
\end{equation*}
Then, after passing to a subsequence,
\begin{equation*}
 \chi_k\overset{*}{\rightharpoonup}\Theta
 \quad\text{in }L^\infty_\loc(D),
 \qquad
 u_k\longrightarrow u
 \quad\text{locally uniformly and strongly in }H^1_\loc(D).
\end{equation*}
The limit belongs to $C^\beta_\loc(D)$ for every $\beta<\alpha$, is harmonic in $\{u\neq0\}$, and satisfies \eqref{eq lambda} with phase $\Theta$.
\end{lemma}

\begin{proof}
Fix $B_{2R}(x_0)\subset\subset D$ and $\beta<\alpha$. For each $k$, \eqref{eq Fbeta} gives \eqref{eq remainderWeiss} with $\gamma=\beta$ at the free boundary points of $u_k$. Notice that this first part requires the initial $C^\alpha$ regularity of each $u_k$ in a qualitative way, but no uniform bound is required. Corollary~\ref{cor weisstoholder}, the uniform $H^1$ bound on $u_k$, and the uniform $L^\infty$ bound in $\chi_k$ give a uniform $C^\beta$ seminorm bound on smaller balls. The mean-value inequality for the subharmonic functions $u_k^2$ combined with the $H^1$ bound also gives a uniform $L^\infty$ bound.

Weak $H^1$ compactness, Arzel\`a--Ascoli, and a diagonal argument now provide a subsequence converging weakly in $H^1_\loc(D)$ and locally uniformly to a function $u$. The uniform convergence shows that $u$ is harmonic in $\{u\neq0\}$. It also follows that $u\in C^\beta_\loc(D)$ for every $\beta<\alpha$. Applying the truncation argument used in the derivation of \eqref{eq outer} to the limit $u$ gives
\begin{equation*}
 \int_D\eta|\nabla u|^2
 =-\int_D u\nabla u\cdot\nabla\eta
 \quad\mbox{for all }\eta\in C_c^\infty(D).
\end{equation*}

For the approximating functions, \eqref{eq outer} and the local uniform convergence imply
\begin{equation*}
 \int_D\eta|\nabla u_k|^2
 =-\int_D u_k\nabla u_k\cdot\nabla\eta
 \longrightarrow-\int_D u\nabla u\cdot\nabla\eta
 =\int_D\eta|\nabla u|^2
 \quad\mbox{for all }\eta\in C_c^\infty(D).
\end{equation*}
Taking nonnegative cutoffs equal to one on a given compact subset and combining this convergence of Dirichlet energies with weak $H^1$ convergence proves that $u_k\to u$ strongly in $H^1_\loc(D)$. Up to extracting a further subsequence, we have that the faces converge in the weak-* sense in $L^\infty$ to a bounded phase $\Theta$. Strong convergence of the gradients lets us pass to the limit in the quadratic terms of \eqref{eq lambda}, while weak-* convergence guarantees the convergence of the phases. Thus $u$ satisfies the same stationarity identity with phase $\Theta$.
\end{proof}

\subsection{Height blow-ups and the frequency bound}\label{subsec frequency}

We first consider $\chi=0$. These are the limits obtained by height normalization when the height grows faster than $r^2$. The argument follows the frequency and dimension-reduction method in \cite[Sections 2.6--2.7]{NRS}; we include the details needed here for signed functions.

\begin{lemma}[Doubling and height monotonicity]\label{lemma doubling}
Suppose that $\chi=0$ and $x_0\in\partial\{u\neq0\}$. Then
\begin{equation}\label{eq doubling}
 H(x_0,r)\leq2^{2N(x_0,r)}H(x_0,r/2).
\end{equation}
If $\gamma=N(x_0,0+)$, then, for $0<s<r$,
\begin{equation}\label{eq heightmonotonicity}
 \frac{H(x_0,s)}{s^{2\gamma}}
 \leq\frac{H(x_0,r)}{r^{2\gamma}}.
\end{equation}
\end{lemma}

\begin{proof}
The first assertion follows by integrating \eqref{eq logheightder}, namely $(\log H(x_0,r))'=2N(x_0,r)/r$, and using the monotonicity of $N(x_0,r)$. For \eqref{eq heightmonotonicity}, we use \eqref{eq logheightder} again to obtain
\begin{equation*}
 \frac{d}{dr}\left(\frac{H(x_0,r)}{r^{2\gamma}}\right)
 =\frac{2H(x_0,r)}{r^{2\gamma+1}}
 \bigl(N(x_0,r)-\gamma\bigr)\geq0.
\end{equation*}
\end{proof}

\begin{lemma}[Tangent functions]\label{lemma prop blowups}
Suppose that $\chi=0$ and let $x_0\in\partial\{u\neq0\}$. For every sequence $r_j\downarrow0$, the height rescalings
\begin{equation}\label{eq HR}
 u_r(x)=\frac{u(x_0+rx)}{\sqrt{H(x_0,r)}}
\end{equation}
have a subsequence converging locally uniformly and strongly in $H^1_\loc(\R^n)$ to a nonzero function $\bar u$. The limit is homogeneous of degree $\gamma=N(x_0,0+)\geq\alpha$, satisfies \eqref{eq lambda} with zero phase, and is harmonic in its nonzero set.

For any nonzero homogeneous solution $v$ of degree $\gamma>0$ in this class and any nonzero free boundary point $z$, we have $N_v(z,0+)\leq N_v(0,0+)=\gamma$. Every tangent of $v$ at $z$ is invariant under translations in direction $z$. Furthermore, if $N_v(z,0+)=\gamma$, then $v$ itself has this invariance.
\end{lemma}

\begin{proof}
Fix $R>0$. For sufficiently small $r$, monotonicity of the frequency and \eqref{eq logheightder} bound
\begin{equation*}
 H_{u_r}(0,R)=\frac{H(x_0,rR)}{H(x_0,r)}
\end{equation*}
uniformly. The same frequency bound controls
$D_{u_r}(0,R)=N_{u_r}(0,R)H_{u_r}(0,R)$, so the rescalings are uniformly bounded in $H^1(B_R)$. Each rescaling satisfies the hypotheses of Lemma~\ref{lemma compactness} with zero phase and the same initial H\"older exponent. That lemma does not require a uniform bound on their initial H\"older seminorms. We can therefore extract a subsequence converging locally uniformly and strongly in $H^1_\loc(\R^n)$ to a function $\bar u$ satisfying the same equations. Since $H_{u_r}(0,1)=1$, local uniform convergence gives $H_{\bar u}(0,1)=1$, so $\bar u$ is nonzero.

For every $\rho>0$, strong convergence of the Dirichlet energies and convergence of the heights give
\begin{equation*}
 N_{\bar u}(0,\rho)
 =\lim_{j\to\infty}N(x_0,\rho r_j)
 =\gamma.
\end{equation*}
The rigidity statement in Lemma~\ref{lemma almgren} shows that $\bar u$ is homogeneous of degree $\gamma$. Moreover, \eqref{eq Fbeta}, with $\chi=0$, gives $N(x_0,r)\geq\beta$ for every $\beta<\alpha$. Letting $r\downarrow0$ and then $\beta\uparrow\alpha$ proves $\gamma\geq\alpha$.

Now let $v$ be a nonzero entire solution homogeneous of degree $\gamma$. Homogeneity and continuity of the frequency with respect to its center imply
\begin{equation*}
 \lim_{R\to\infty}N_v(z,R)
 =\lim_{R\to\infty}N_v(z/R,1)
 =N_v(0,1)=\gamma.
\end{equation*}
Thus monotonicity gives
\begin{equation}\label{eq monfr}
 N_v(z,0+)\leq N_v(z,R)\leq\gamma
 \qquad\mbox{for every }R>0.
\end{equation}

Consider a height tangent of $v$ at $z$, obtained from the rescalings
$v_r(x)=v(z+rx)/\sqrt{H_v(z,r)}$. Euler's identity for $v$ yields
\begin{equation*}
 z\cdot\nabla v_r
 =r\bigl(\gamma v_r-x\cdot\nabla v_r\bigr).
\end{equation*}
The rescalings are locally bounded in $H^1$ by the first part of the proof, so the right hand side tends to zero in $L^2_\loc$ as $r\downarrow0$. Every tangent is therefore invariant under translations in direction $z$.

Finally, if $N_v(z,0+)=\gamma$, \eqref{eq monfr} implies that $N_v(z,R)=\gamma$ for every $R>0$. Hence $v$ is homogeneous of degree $\gamma$ about both $0$ and $z$. Subtracting the identities
$x\cdot\nabla v=\gamma v$ and $(x-z)\cdot\nabla v=\gamma v$ gives $z\cdot\nabla v=0$, which implies the claimed invariance.
\end{proof}

\begin{lemma}[Lower frequency bound when $\chi=0$]\label{lemma dir}
If $\chi=0$, then at every free boundary point $z\in\partial\{u\neq0\}$,
\begin{equation}\label{eq lowerfreqbound}
 N(z,0+)\geq1.
\end{equation}
Moreover, if a nonzero global solution is homogeneous of degree one, then
$|u(x)|=c|x\cdot e|$ for some $c>0$ and some unit vector $e$.
\end{lemma}

\begin{proof}
By Lemma~\ref{lemma prop blowups}, a height blow-up at any free boundary point is nonzero and homogeneous, with degree equal to the limiting frequency there. It therefore suffices to prove the frequency bound for nonzero homogeneous solutions. We prove this bound and the rigidity statement together by induction on the dimension.

In dimension two, we can apply the planar argument of \cite[Lemma 2.20]{NRS} to $|v|$. Indeed, $|v|$ is harmonic in its positivity set, and its stress tensor agrees almost everywhere with that of $v$. The angular equation in each component of $\{v\neq0\}$, together with stationarity across the intervening rays, gives, up to rotation,
\begin{equation}\label{eq 2d}
 |v(r,\theta)|
 =c r^{N/2}\left|\sin\left(\frac{N\theta}{2}\right)\right|,
 \qquad N\geq2,
\end{equation}
for some $c>0$ and integer $N$. In particular, $\gamma=N/2\geq1$. Equality holds precisely when $|v(x)|=c|x\cdot e|$.

Suppose that both assertions hold in dimension $n-1$, with $n\geq3$. Let $v$ be a nonzero homogeneous solution in $\R^n$ of degree $0<\gamma<1$. There must be a nonzero point $z\in\partial\{v\neq0\}$. Otherwise the origin would be its only zero: $v$ would have one sign and be harmonic away from the origin. Since it is continuous there, the singularity would be removable, and the minimum principle would contradict $v(0)=0$.

Take a height tangent of $v$ at such a point $z$. By Lemma~\ref{lemma prop blowups}, it is nonzero, homogeneous of degree $N_v(z,0+)\leq\gamma<1$, and invariant in direction $z$. It therefore defines a solution in $\R^{n-1}$ of the same degree. The reduced stationarity identity follows by testing with vector fields given by a field in $\R^{n-1}$ times a cutoff in the invariant variable. Harmonicity outside the zero set and local H\"older continuity pass to the reduced function as well. This contradicts the induction hypothesis and proves \eqref{eq lowerfreqbound} in dimension $n$.

Finally, let $v$ be homogeneous of degree one. As above, it has a nonzero free boundary point $z$. The bound just proved and Lemma~\ref{lemma prop blowups} give
\begin{equation*}
 1\leq N_v(z,0+)\leq1.
\end{equation*}
The equality case of that lemma shows that $v$ itself is invariant in direction $z$. Applying the rigidity statement in dimension $n-1$ to the reduced solution completes the induction.
\end{proof}

\subsection{H\"older improvement}\label{subsec improvement}

Armed with the results proved in this section, we show that $F(z,r)\geq1$ for $z\in\partial\{u\neq0\}$, where $F$ is the generalized frequency defined in \eqref{eq newmon}. This observation will allow us to improve the initial $\alpha$-H\"older regularity to $\gamma$-H\"older regularity for every $\gamma\in(0,1)$.

\begin{lemma}\label{lemma Fone}
At every free boundary point $z$ and every admissible radius,
\begin{equation}\label{eq Fone}
 F(z,r)\geq1.
\end{equation}
\end{lemma}

\begin{proof}
We argue by contradiction and assume that $F(z,r_0)=q<1$ for some admissible radius $r_0$. By Theorem~\ref{thm monF}, $F(z,r)$ is nondecreasing on $(0,r_0)$ and converges as $r\downarrow0$ to some $\theta\leq q$. Moreover, the initial $\alpha$-H\"older regularity of $u$, combined with \eqref{eq Fbeta}, gives $F(z,r)\geq\beta$ for every $\beta<\alpha$, so $\theta\geq\alpha$.

Since $N(z,r)\leq F(z,r)\leq q$ on $(0,r_0)$, \eqref{eq logheightder} gives the differential inequality
\begin{equation*}
 \frac{d}{dr}H(z,r)\leq\frac{2q}{r}H(z,r).
\end{equation*}
Integrating between $0<r<s\leq r_0$ gives
\begin{equation}\label{eq heightlb}
 H(z,r)\geq H(z,s)\left(\frac{r}{s}\right)^{2q}.
\end{equation}
Consequently, setting $s=r_0$, we find that $H(z,r)\geq H(z,r_0)(r/r_0)^{2q}$, so that
\begin{equation*}
 0\leq F(z,r)-N(z,r)
 =\frac{r^{2-n}\int_{B_r(z)}\chi}{H(z,r)}
 =O(r^{2-2q})\longrightarrow0,
\end{equation*}
since $q<1$. We conclude that
\begin{equation}\label{eq convN}
 \lim_{r\downarrow0}N(z,r)=\theta.
\end{equation}

Consider now the height rescalings $ u_r(x)=\frac{u(z+rx)}{\sqrt{H(z,r)}}$. They satisfy \eqref{eq lambda} with phases
$\frac{r^2}{H(z,r)}\chi(z+r\,\cdot)$, which converge to zero locally in $L^\infty$ by \eqref{eq heightlb}. On the other hand, \eqref{eq convN} and the height normalization $H_{u_r}(0,1)=1$ provide uniform local $H^1$ bounds: \eqref{eq logheightder} controls the height at each fixed radius, and $D_{u_r}=N_{u_r}H_{u_r}$ controls the Dirichlet energy. By Lemma~\ref{lemma compactness}, we can find a subsequence converging locally uniformly and strongly in $H^1_\loc$ to a function $v$ with $\chi=0$.

The height normalization combined with \eqref{eq heightlb} implies that $v$ is nontrivial. Strong convergence and \eqref{eq convN} therefore give
\begin{equation*}
 N_v(0,\rho)=\theta
 \qquad\mbox{for every }\rho>0.
\end{equation*}
Finally, Lemma~\ref{lemma dir} implies $\theta\geq1$, contradicting $\theta\leq q<1$.
\end{proof}

We finish this section with the proof of Theorem~\ref{thm anyholder}.

\begin{proof}[Proof of Theorem~\ref{thm anyholder}]
The lower bound on the generalized frequency from
Lemma~\ref{lemma Fone} implies that $W(z,r)\geq0$ at every free boundary point. In particular, the hypothesis of Corollary~\ref{cor weisstoholder} holds for every $\gamma<1$, and hence $u\in C^\gamma_\loc(\Om)$ for every such $\gamma$.

We now use the bound $W\geq0$ to obtain the finer modulus stated in the theorem. Differentiating the height gives
\begin{equation}\label{eq derfrequency}
 \frac{d}{dr}\left(\frac{H(z,r)}{r^2}\right)
 =\frac{2}{r^3}\bigl(D(z,r)-H(z,r)\bigr)
 =\frac2r\left(W(z,r)-r^{-n}\int_{B_r(z)}\chi\right).
\end{equation}
Therefore, for $0<r<r_0$,
\begin{equation}\label{eq rlogr}
 \frac{H(z,r)}{r^2}
 \leq\frac{H(z,r_0)}{r_0^2}
 +2|B_1|\|\chi\|_{L^\infty(B_{r_0}(z))}
 \log\frac{r_0}{r}.
\end{equation}
The logarithm appears in the estimate for the height, which controls the square of the size of $u$. Taking its square root gives the $\sqrt{\log}$ factor.

We now pass from \eqref{eq rlogr}, centered at free boundary points, to an estimate for arbitrary pairs of points, as in Lemma~\ref{lemma reg criterion}. Fix $B_{2R}(x_0)\subset\subset\Om$. For free boundary points $z\in B_{3R/2}(x_0)$, take $r_0=R/4$. The mean-value inequality for the subharmonic function $u^2$ and \eqref{eq rlogr} give, for $0<r<R/8$,
\begin{equation}\label{eq growthlog}
 \sup_{B_{r/2}(z)}|u|
 \leq Cr\left(
 \frac{\|u\|_{L^\infty(B_{2R}(x_0))}}{R}
 +\|\chi\|_{L^\infty(B_{2R}(x_0))}^{1/2}
 \sqrt{1+\log\frac{R}{r}}
 \right).
\end{equation}

Let $x\in B_R(x_0)\cap\{u\neq0\}$ be at distance $0<d<R/64$ from the zero set. A nearest zero belongs to the boundary of the phase containing $x$. Applying \eqref{eq growthlog} at this free boundary point with a radius comparable to $d$, and then using an interior estimate for the harmonic function $u$ in $B_d(x)$, we obtain
\begin{equation}\label{eq gradientlog}
 |\nabla u(x)|
 \leq C\left(
 \frac{\|u\|_{L^\infty(B_{2R}(x_0))}}{R}
 +\|\chi\|_{L^\infty(B_{2R}(x_0))}^{1/2}
 \sqrt{1+\log\frac{R}{d}}
 \right).
\end{equation}
If the distance from $x$ to the zero set is at least $R/64$, the usual interior harmonic estimate gives the corresponding bound with $R$ in place of $d$.

Finally, let $x,y\in B_R(x_0)$ and put $s=|x-y|$. If $s$ is comparable to $R$, the supremum bound for $u$ suffices. Suppose that $s$ is smaller. If one endpoint lies within distance $2s$ of the zero set, then both endpoints lie within distance $3s$ of it. Applying \eqref{eq growthlog} at nearby free boundary points bounds both values by the right hand side of \eqref{eq anyHolderbound}. Otherwise, the segment joining $x$ and $y$ stays at distance at least $s$ from the zero set. We integrate \eqref{eq gradientlog} along this segment and obtain the same estimate. This proves \eqref{eq anyHolderbound}.
\end{proof}
\section{Lipschitz regularity}\label{sec lipschitz}

For the rest of this section, we assume that $u$ and $\chi$ satisfy the hypotheses of Theorem~\ref{thm main}. In particular, $u\geq0$ and
$\one_{\{u>0\}}\leq\chi\leq1$. We write $\omega_n=|B_1|$ and $W=W_1$. Our ultimate goal is to obtain a uniform local bound on $H(z,r)/r^2$, which would give the desired Lipschitz estimate by Lemma~\ref{lemma reg criterion}. We first show, through a chain of auxiliary results, that this quotient remains bounded as $r\downarrow0$ at each fixed free boundary point.

\begin{lemma}\label{lemma fulldens}
Fix $z\in\partial\{u>0\}$ and $r_0>0$ such that
$B_{r_0}(z)\subset\subset\Om$. Given $\varepsilon>0$, there exists $A>0$, depending on $u,z$, and $r_0$, such that, if $0<r<r_0/2$ and
\begin{equation*}
 \frac{H(z,r)}{r^2}>A,
\end{equation*}
then
\begin{equation*}
 \omega_n-\frac1{r^n}\int_{B_r(z)}\chi<\varepsilon.
\end{equation*}
In fact, along any sequence for which $H(z,r_k)/r_k^2\to\infty$, the density of $\{u>0\}$ in $B_{r_k}(z)$ tends to one.
\end{lemma}

\begin{proof}
We argue by contradiction and consider a sequence $r_k$ such that
$H(z,r_k)/r_k^2\to\infty$ but
\begin{equation*}
 \omega_n-\frac1{r_k^n}\int_{B_{r_k}(z)}\chi\geq\varepsilon_0
\end{equation*}
for some $\varepsilon_0>0$. After passing to a subsequence, necessarily $r_k\downarrow0$. By Lemma~\ref{lemma Fone} and Weiss monotonicity, we have $0\leq W(z,r)\leq W(z,r_0)$. Combining this with $0\leq\chi\leq1$ and \eqref{eq derfrequency} gives
\begin{equation}\label{eq derfrequencybound}
 \left|\frac{d}{dr}\left(\frac{H(z,r)}{r^2}\right)\right|
 \leq\frac{2\bigl(W(z,r_0)+\omega_n\bigr)}r
\end{equation}
for almost every $r\in(0,r_0)$. Integrating between $r_k$ and $tr_k$, with $t>0$ fixed, we deduce that
\begin{equation}\label{eq slogs}
 \frac{H(z,tr_k)}{H(z,r_k)}\longrightarrow t^2.
\end{equation}

Now consider the height rescalings
$v_k(x)=u(z+r_kx)/\sqrt{H(z,r_k)}$, whose corresponding phases
$r_k^2\chi(z+r_k\,\cdot)/H(z,r_k)$ tend to zero locally in $L^\infty$. Moreover, we have the identity
\begin{equation}\label{eq idfreq}
 \frac{D(z,r)}{H(z,r)}-1
 =\frac{r^2}{H(z,r)}
 \left(W(z,r)-r^{-n}\int_{B_r(z)}\chi\right).
\end{equation}
The boundedness of $W$ and \eqref{eq slogs} show that
$N(z,tr_k)\to1$ for every fixed $t>0$. Thus the Dirichlet energies of the rescalings are locally bounded; in fact,
\begin{equation*}
 D_{v_k}(0,t)=N(z,tr_k)H_{v_k}(0,t)\longrightarrow t^2.
\end{equation*}
Lemma~\ref{lemma compactness} gives, after passing to a subsequence, a locally uniform and strong $H^1_\loc$ limit $v$ with zero phase. By our normalization and \eqref{eq slogs}, we have $H_v(0,t)=t^2$. The convergence of the Dirichlet energies gives $D_v(0,t)=t^2$. Thus $N_v(0,t)=1$ for every $t>0$. By Lemmas~\ref{lemma almgren} and \ref{lemma dir}, and since $v\geq0$, we have $v(x)=c|x\cdot e|$ for some $c>0$ and some unit vector $e$.

In particular, $v>0$ almost everywhere. Local uniform convergence implies
\begin{equation*}
 |\{v_k>0\}\cap B_1|\longrightarrow\omega_n.
\end{equation*}
Since $\one_{\{u>0\}}\leq\chi\leq1$, we have
\begin{equation*}
 |\{v_k>0\}\cap B_1|
 \leq\frac1{r_k^n}\int_{B_{r_k}(z)}\chi
 \leq\omega_n.
\end{equation*}
The phase averages therefore tend to $\omega_n$, giving the desired contradiction.
\end{proof}

We now prove the quadratic bound for the height function.

\begin{lemma}\label{lemma finheight}
For every $z\in\partial\{u>0\}$,
\begin{equation*}
 \limsup_{r\downarrow0}\frac{H(z,r)}{r^2}<\infty.
\end{equation*}
\end{lemma}

\begin{proof}
\medskip
\noindent\emph{Preparations.}
We start by deriving a differential inequality that will be used in the last part of the proof. Direct differentiation gives the following alternative expression for \eqref{eq derfrequency}:
\begin{equation}\label{eq derfrequency1}
 \frac{d}{dr}\left(\frac{H(z,r)}{r^2}\right)
 =\frac{2}{r^{n+1}}
 \int_{\partial B_r(z)}u\left(\partial_\nu u-\frac ur\right).
\end{equation}
On the other hand, Weiss monotonicity gives
\begin{equation}\label{eq derweiss}
 \frac{d}{dr}W(z,r)
 =\frac2{r^n}\int_{\partial B_r(z)}
 \left(\partial_\nu u-\frac ur\right)^2.
\end{equation}
Combining \eqref{eq derfrequency1}, \eqref{eq derweiss}, and Cauchy--Schwarz yields
\begin{equation}\label{eq diffineqH}
 \left(\frac{d}{dr}\left(\frac{H(z,r)}{r^2}\right)\right)^2
 \leq\frac{2H(z,r)}{r^3}\frac{d}{dr}W(z,r).
\end{equation}
As noted in the previous proof, Lemma~\ref{lemma Fone} and Weiss monotonicity imply that the limit $W(z,0+)$ exists and belongs to $[0,\infty)$. We consider two cases.

\medskip
\noindent\emph{Case 1: $W(z,0+)\geq\omega_n$.}
Monotonicity gives $W(z,r)\geq\omega_n$, whereas
$r^{-n}\int_{B_r(z)}\chi\leq\omega_n$. Thus \eqref{eq derfrequency} shows that $r\mapsto H(z,r)/r^2$ is nondecreasing. In particular, it remains bounded as $r\downarrow0$.

\medskip
\noindent\emph{Case 2: $W(z,0+)<\omega_n$.}
We argue by contradiction and assume that $H(z,r_k)/r_k^2\to\infty$ for some sequence $r_k\downarrow0$. By Lemma~\ref{lemma fulldens},
\begin{equation}\label{eq seqdens}
 \lim_{k\to\infty}\frac1{r_k^n}\int_{B_{r_k}(z)}\chi
 =\omega_n.
\end{equation}
We claim, more precisely, that
\begin{equation}\label{eq infH}
 \lim_{r\downarrow0}\frac{H(z,r)}{r^2}=\infty
\end{equation}
and
\begin{equation}\label{eq fullchi}
 \lim_{r\downarrow0}\frac1{r^n}\int_{B_r(z)}\chi=\omega_n.
\end{equation}

To see this, put $\delta=(\omega_n-W(z,0+))/4$ and take the threshold $A$ given by Lemma~\ref{lemma fulldens} with $\varepsilon=\delta$. We choose $k$ large enough that $H(z,r_k)/r_k^2>A$ and
\begin{equation*}
 W(z,r)\leq W(z,0+)+\delta
 \qquad\mbox{for }0<r\leq r_k.
\end{equation*}
At every radius in this interval for which $H(z,r)/r^2>A$, Lemma~\ref{lemma fulldens} gives
\begin{equation*}
 \frac1{r^n}\int_{B_r(z)}\chi>\omega_n-\delta.
\end{equation*}
Consequently, \eqref{eq derfrequency} implies
\begin{equation}\label{eq derr0}
 \frac{d}{dr}\left(\frac{H(z,r)}{r^2}\right)
 \leq-\frac{4\delta}{r}
\end{equation}
for almost every such radius.

We observe now that the previous inequality propagates the condition $H(z,r)/r^2>A$ towards smaller radii. Indeed, on the interval where this condition holds and which contains $r_k$, integration of \eqref{eq derr0} towards smaller radii makes the quotient strictly larger. Hence the interval cannot have a positive left endpoint at which the quotient equals $A$. It follows that \eqref{eq derr0} holds almost everywhere on $(0,r_k)$. Integrating it proves \eqref{eq infH}, while Lemma~\ref{lemma fulldens} then gives \eqref{eq fullchi}.

Finally, \eqref{eq diffineqH} and \eqref{eq derr0} imply
\begin{equation}\label{eq diffineqW}
 \frac{d}{dr}W(z,r)
 \geq\frac{8\delta^2}{r\bigl(H(z,r)/r^2\bigr)}
 \qquad\mbox{for almost every }r\in(0,r_k).
\end{equation}
On the other hand, \eqref{eq rlogr} gives
\begin{equation*}
 \frac{H(z,r)}{r^2}
 \leq\frac{H(z,r_k)}{r_k^2}
 +2\omega_n\log\frac{r_k}{r}.
\end{equation*}
Combining the last two inequalities, we deduce
\begin{equation*}
 \frac{d}{dr}W(z,r)
 \geq
 \frac{8\delta^2}
 {r\left(H(z,r_k)/r_k^2+2\omega_n\log(r_k/r)\right)}.
\end{equation*}
The right hand side is not integrable at zero, contradicting the finiteness of $W(z,r_k)-W(z,0+)$.
\end{proof}

With the pointwise bound on $H(z,r)/r^2$, we can show that the generalized phase density exists at every free boundary point. The argument is close to \cite[Proposition 3.5]{KW}.

\begin{lemma}[Existence of generalized density]\label{lemma density}
At every $z\in\partial\{u>0\}$,
\begin{equation}\label{eq density}
 \lim_{r\downarrow0}r^{-n}\int_{B_r(z)}\chi
 =W(z,0+)\in[0,\omega_n].
\end{equation}
Moreover, every sequence of linear rescalings $u(z+r_j\,\cdot)/r_j$, with $r_j\downarrow0$, has a subsequence converging locally uniformly and strongly in $H^1_\loc(\R^n)$ to a one-homogeneous function, possibly zero.
\end{lemma}

\begin{proof}
By Lemma~\ref{lemma finheight}, the heights of the linear rescalings
$u_r(x)=u(z+rx)/r$ are bounded on each fixed ball. Subharmonicity gives the corresponding local $L^2$ bounds. Since $W(z,r)$ is nondecreasing and bounded above at a fixed radius, its definition also bounds the Dirichlet energies of $u_r$. We may therefore apply Lemma~\ref{lemma compactness}. Given any sequence $r_j\downarrow0$, a subsequence converges locally uniformly and strongly in $H^1_\loc$ to a function $v$. Along a further subsequence, the phases $\chi(z+r_j\,\cdot)$ converge weakly-* to a phase $\Theta$ satisfying $\one_{\{v>0\}}\leq\Theta\leq1$.

We next show that $v$ is one-homogeneous. For $0<a<b$, we integrate \eqref{eq weissbeta derivative} with $\beta=1$ to find
\begin{equation*}
 W(z,br_j)-W(z,ar_j)
 =2\int_{B_b\setminus B_a}|x|^{-n}
 \left(\partial_r u_{r_j}-\frac{u_{r_j}}{|x|}\right)^2dx
 \longrightarrow0.
\end{equation*}
Strong $H^1$ convergence shows that the corresponding integral vanishes for $v$. Since $a$ and $b$ are arbitrary, $v$ is one-homogeneous; it may be identically zero. Also, $v\Delta v=0$ locally in the sense of distributions. We can therefore combine one-homogeneity with \eqref{eq testdir}, applied to $v$, to deduce
\begin{equation*}
 \int_{B_1}|\nabla v|^2=\int_{\partial B_1}v^2.
\end{equation*}
Passing to the limit in the definition of $W(z,r_j)$, we conclude that
\begin{equation*}
 W(z,0+)
 =\lim_{j\to\infty}r_j^{-n}\int_{B_{r_j}(z)}\chi
\end{equation*}
along the extracted subsequence. Every sequence has a subsequence with this same limit, proving \eqref{eq density}. The range of the limit follows from $0\leq\chi\leq1$.
\end{proof}

With the existence of the generalized density at every free boundary point, we can refine the dichotomy used in Lemma~\ref{lemma finheight}. At points where the phase has maximal density, the height satisfies a monotonicity formula.

\begin{lemma}[Height monotonicity at full phase density]\label{lemma fullphase}
Let $z\in\partial\{u>0\}$ and suppose that
\begin{equation*}
 \lim_{r\downarrow0}r^{-n}\int_{B_r(z)}\chi=\omega_n.
\end{equation*}
Then $r\mapsto H(z,r)/r^2$ is nondecreasing on its interval of definition.
\end{lemma}

\begin{proof}
By Lemma~\ref{lemma density} and Weiss monotonicity,
$W(z,r)\geq\omega_n$. Since
$r^{-n}\int_{B_r(z)}\chi\leq\omega_n$, the conclusion follows from \eqref{eq derfrequency}.
\end{proof}

At points where the generalized density is less than $\omega_n$, we obtain a complementary bound using a viscosity argument. We show that a smooth function touching $u$ from below at such a point cannot have slope greater than one.

\begin{lemma}[Viscosity supersolution at points of nonmaximal density]\label{lemma touching}
Let $z\in\partial\{u>0\}$ and suppose that
\begin{equation*}
 \lim_{r\downarrow0}r^{-n}\int_{B_r(z)}\chi<\omega_n.
\end{equation*}
If a smooth function $\phi$ touches $u$ from below at $z$, then
$|\nabla\phi(z)|\leq1$.
\end{lemma}

\begin{proof}
Since $\phi$ touches $u$ from below at $z$, we have $\phi\leq u$ near $z$ and $\phi(z)=u(z)=0$. Consider the linear rescalings of $\phi$ and $u$ at $z$. By Lemma~\ref{lemma density}, a subsequence of the rescalings of $u$ converges locally uniformly and strongly in $H^1_\loc$ to a one-homogeneous function $v$. If $\nabla\phi(z)=0$, there is nothing to prove. Otherwise $v$ is nonzero and $v(x)\geq\nabla\phi(z)\cdot x$ for $x\in\R^n$. In particular, $v>0$ in the half-space $\{x\cdot e>0\}$, where $e=\nabla\phi(z)/|\nabla\phi(z)|$. By homogeneity, $v(r,\theta)=r\psi(\theta)$, with $\psi$ satisfying $-\Delta_{\mathbb S^{n-1}}\psi=(n-1)\psi$ in each connected component of $\{\psi>0\}$ on $\mathbb S^{n-1}$. By strict monotonicity of the first Dirichlet eigenvalue, the component containing the hemisphere $\{x\cdot e>0\}$ must coincide with this hemisphere. It follows that $v=c(x\cdot e)$ in that half-space for some $c\geq|\nabla\phi(z)|$. Any other positivity component of $\{\psi>0\}$ would lie in the opposite hemisphere. Since its first eigenvalue is also $n-1$, the same comparison forces it to fill that hemisphere. However, this would imply that $v>0$ almost everywhere. Since $\Theta=1$ wherever $v>0$, we would obtain
$\lim_{r\downarrow0}r^{-n}\int_{B_r(z)}\chi=\omega_n$, contradicting the hypothesis. Thus $v(x)=c(x\cdot e)_+$.

Along the chosen subsequence, the rescaled phases $\chi(z+r_j\,\cdot)$ converge weakly-* to $\Theta$, and $v$ satisfies \eqref{eq stationarity} with phase $\Theta$. We have $\Theta=1$ in $\{x\cdot e>0\}$. In the opposite half-space, $\nabla v=0$, so \eqref{eq stationarity} gives
\begin{equation*}
 \int_{\{x\cdot e<0\}}\Theta\,\diver T=0
 \quad\mbox{for all }T\in\Lip_c(\{x\cdot e<0\};\R^n).
\end{equation*}
Hence $\Theta$ equals a constant $\theta\in[0,1]$ there. Testing the inner variation with vector fields supported in balls intersecting $\{x\cdot e=0\}$ gives the transmission condition
\begin{equation*}
 1-c^2=\theta.
\end{equation*}
Therefore $|\nabla\phi(z)|\leq c\leq1$, as claimed.
\end{proof}

We combine the last two lemmas to derive the Lipschitz estimate. At points of full generalized density we use height monotonicity; at the remaining points we use the touching condition and a standard annular barrier argument, as in \cite[Lemma 1.5]{CafarelliSalsa}.

\begin{proof}[Proof of Theorem~\ref{thm main}]
Fix $B_{2R}(x_0)\subset\subset\Om$ and $y\in B_R(x_0)\cap\{u>0\}$. If the distance from $y$ to the zero set is at least $R/8$, the interior gradient estimate gives
\begin{equation*}
 |\nabla u(y)|
 \leq\frac{C}{R}\|u\|_{L^\infty(B_{2R}(x_0))}.
\end{equation*}
Otherwise, write $d=\operatorname{dist}(y,\{u=0\})<R/8$ and choose $z\in\partial\{u>0\}$ with $|y-z|=d$. Lemma~\ref{lemma density} gives the existence of the generalized density at $z$. We consider two cases.

\medskip
\noindent\emph{Case 1: The generalized density at $z$ is $\omega_n$.}
Combining Lemma~\ref{lemma fullphase} with the mean-value inequality for the subharmonic function $u^2$, we obtain
\begin{equation*}
 u(y)^2
 \leq C H(z,4d)
 \leq C\left(\frac{d}{R}\right)^2H(z,R/2)
 \leq C\frac{d^2}{R^2}
       \|u\|_{L^\infty(B_{2R}(x_0))}^2.
\end{equation*}

\medskip
\noindent\emph{Case 2: The generalized density at $z$ is less than $\omega_n$.}
Set $s=d/2$. By Harnack's inequality in $B_{2s}(y)$, there is a dimensional constant $C$ such that
\begin{equation*}
 u\geq m_s:=\frac1C\|u\|_{L^\infty(B_s(y))}
 \qquad\mbox{in }B_s(y).
\end{equation*}
Let $\phi$ be the radial harmonic function in
$B_2\setminus\overline{B_1}$ which equals one on $\partial B_1$ and zero on $\partial B_2$. The function
\begin{equation*}
 \psi_s(x)=m_s\phi\left(\frac{x-y}{s}\right)
\end{equation*}
is bounded above by $u$ on the boundary of
$B_{2s}(y)\setminus\overline{B_s(y)}$, and hence throughout the annulus. The radial formula for $\phi$ is an affine combination of a fundamental solution of the Laplacian and a constant. It extends smoothly across the outer sphere and is negative just outside it. Thus $\psi_s$ is defined in a neighborhood of $z$ and touches $u$ from below there. Since
\begin{equation*}
 |\nabla\psi_s(z)|=c(n)\frac{m_s}{s},
\end{equation*}
Lemma~\ref{lemma touching} gives $m_s\leq Cs$. In particular,
\begin{equation*}
 u(y)\leq\|u\|_{L^\infty(B_s(y))}\leq Cs.
\end{equation*}

In either case,
\begin{equation*}
 u(y)\leq Cd\left(
 1+\frac{\|u\|_{L^\infty(B_{2R}(x_0))}}R
 \right).
\end{equation*}
Harnack's inequality and the interior gradient estimate in $B_d(y)$ give
\begin{equation*}
 |\nabla u(y)|
 \leq\frac{C}{d}\sup_{B_{d/2}(y)}u
 \leq\frac{C}{d}u(y)
 \leq C\left(
 1+\frac{\|u\|_{L^\infty(B_{2R}(x_0))}}R
 \right).
\end{equation*}
Together with the estimate at points whose distance from the zero set
is at least $R/8$, this proves the desired gradient bound throughout
$\{u>0\}\cap B_R(x_0)$. 
\end{proof}

\end{document}